\baselineskip=1.5\baselineskip

\magnification=\magstep1
\catcode`@=11
\font\teneufm=eufm10
\font\seveneufm=eufm7
\font\fiveeufm=eufm5
\newfam\eufmfam
\textfont\eufmfam=\teneufm
\scriptfont\eufmfam=\seveneufm
\scriptscriptfont\eufmfam=\fiveeufm
\def\frak{%
 \let\next\relax
 \ifmmode\let\next\frak@
 \else\def\next{\errmessage{Use \string\frak\space in math mode only}}
 \fi
 \next}
\def\frak@#1{{\frak@@{#1}}}
\def\frak@@#1{\fam\eufmfam#1}
\let\goth\frak
\catcode`\@=12

\def\bbbr{{\rm I\!R}}

\def\bbbn{{\rm I\!N}}

\def\bbbzz{{\rm Z\!\!\!Z}}
\def\bbbc{{\mathchoice {\setbox0=\hbox{$\displaystyle\rm
C$}\hbox{\hbox
to0pt{\kern0.4\wd0\vrule height0.9\ht0\hss}\box0}}
{\setbox0=\hbox{$\textstyle\rm C$}\hbox{\hbox
to0pt{\kern0.4\wd0\vrule height0.9\ht0\hss}\box0}}
{\setbox0=\hbox{$\scriptstyle\rm C$}\hbox{\hbox
to0pt{\kern0.4\wd0\vrule height0.9\ht0\hss}\box0}}
{\setbox0=\hbox{$\scriptscriptstyle\rm C$}\hbox{\hbox
to0pt{\kern0.4\wd0\vrule height0.9\ht0\hss}\box0}}}}

\def\nz{\hfill\break\noindent}
\def\sn{\smallskip\noindent}
\def\mn{\medskip\noindent}
\def\bn{\bigskip\noindent}
\font\titefont=cmbx9 scaled \magstep2

\def\qed{\vrule height 1.2ex width 1.1ex depth -.1ex}
\def\litem{\par\noindent
               \hangindent=\parindent\ltextindent}
\def\ltextindent#1{\hbox to\hangindent{#1\hss}\ignorespaces}

\def\={\!=\!}
\def\A{{\cal{A}}}
\def\M{{\cal{M}}}
\def\N{{\cal{N}}}
\def\H{{\cal{H}}}
\def\D{{\cal{D}}}
\def\gf{{\goth{f}}}
\def\cP{{\cal{P}}}
\def\cQ{{\cal{Q}}}
\def\gg{{\goth{g}}}
\def\h{{\goth{h}}}
\def\k{{\goth{k}}}

\def\gK{{\goth{K}}}
\def\X{{\cal{X}}}
\def\gH{{\goth{H}}}
\def\gK{{\goth{K}}}

\def\gN{{\goth{N}}}

\noindent
\vskip1cm \noindent 
\centerline{\titefont{Operator Representations of a ${\bf q}$-Deformed 
Heisenberg Algebra}}
\mn
\centerline{{Konrad Schm{\"u}dgen}}
\centerline{Universit{\"a}t Leipzig,}
\centerline{Fakult{\"a}t f{\"u}r Mathematik und Informatik und NTZ,}
\centerline{Augustusplatz 10/11, D--04109 Leipzig, Germany}
\bn
{\bf Abstract.}~ {\it A class of well-behaved $\ast$-representations of a $q$-deformed Heisenberg
algebra introduced in refs. 10 and 3 is studied and classified.}\sn

The idea to develop a $q$-deformed quantum mechanics by using quantum groups
has been investigated in several papers$^{2,3,6,11,13}$. Such approaches are usually based on a $q$-deformed phase space algebra
which is derived from the noncommutative differential calculus of the 
$q$-deformed configuration space$^{7,14}$. Following the standard procedure in 
quantum mechanics one has to represent the $q$-deformed position and
momentum operators by essentially self-adjoint operators acting on a Hilbert
space. More precisely, one has to find appropriate $\ast$-representations of
the phase space $\ast$-algebra by unbounded operators in a Hilbert space.
In the case of general Euclidean or Minikowski phase spaces the study and
classification of these $\ast$-representations turns out to be technically 
complicated because of the many relations and also because of the various
difficulties concerning unbounded operators.

The aim of this paper is to give a rigorous treatment of well-behaved operator 
representations for one of the simplest example - the one-dimensional 
$q$-deformed Heisenberg algebra which was invented in refs. 11 and 3. 
Representations of this algebra have been investigated in ref. 3. 
Since this $\ast$-algebra occurs as a subalgebra of other larger $\ast$-algebras, 
the study of general not necessarily irreducible $\ast$-representations
seems to be important as well. We shall develop and analyze an operator-theoretic
model for such general representations of the $q$-deformed Heisenberg
algebra. This model might be used as a tool kit for the study of 
representations of larger $\ast$-algebras.\sn 

This paper is organized as follows. Section I contains the definition and some simple 
algebraic properties of the $q$-deformed Heisenberg algebra $\A(q)$. In
Section II we develop a general operator-theoretic model for certain triples of operators
which will lead in Section V to representations of the $\ast$-algebra 
$\A(q)$. In Section III the irreducibility and the unitary equivalence of these
operator triples are investigated and a number of examples are treated.
In Section IV we give a characterization of these operator triples
by a number of natural conditions. In Section V we define the
self-adjoint $\ast$-representations of the $\ast$-algebra $\A(q)$ obtained
by means of these operator triples.

In a forthcoming paper we shall study the spectrum of the operator $X$.
For this analysis the $q$-Fourier transform$^{5,4}$ will play a crucial role.
\bn
{\bf I. The ${\bf q}$-Heisenberg algebra}
\sn
For a positive real number $q\ne 1$, let $\A(q)$ denote the complex
unital algebra with four generators ${\bf p,x,u,u^{-1}}$ subject to the defining 
relations
$$
{\bf up}=q~{\bf pu},~{\bf ux}=q^{-1} {\bf xu},~{\bf uu}^{-1}={\bf u}^{-1} 
{\bf u}=1,\eqno(1)
$$
$$
{\bf px}-q~{\bf xp}={\rm i}(q^{3/2}-q^{-1/2}) {\bf u},~{\bf xp}-q~ {\bf px}=
-{\rm i}~(q^{3/2}-q^{-1/2}) {\bf u}^{-1}~~,\eqno(2)
$$
where ${\rm i}$ denotes the imaginary unit. An equivalent set of relations is obtained if (2) is replaced by
$$
{\bf px}={\rm i}~ q^{1/2} {\bf u}^{-1} - {\rm i}~ q^{-1/2} {\bf u},~{\bf xp}=
{\rm i}~ q^{-1/2} {\bf u}^{-1}-{\rm i}~ q^{1/2} {\bf u}.\eqno (2)^\prime
$$
From (1) and (2)$^\prime$ it follows that the set of 
elements $\{{\bf p^ru^n,x^su^n}; r\in\bbbn_0, s\in\bbbn, n\in\bbbzz\}$ 
is a vector space basis of $\A(q)$. 

The algebra $\A(q)$ becomes a $\ast$-algebra with involution defined
on the generators by 
$$
{\bf p}={\bf p^\ast},~{\bf x}={\bf x^\ast},~{\bf u^\ast}={\bf u^{-1}} ~.\eqno(3)
$$
Indeed, it suffices to check that the defining relations (1) and (2)$^\prime$
of $\A(q)$ are invariant under the involution (3) which is easily done.

From (1), (2)$^\prime$ and (3) we conclude that there are $\ast$-isomorphisms
$\rho_1$ and $\rho_2$ of the $\ast$-algebras $\A(q)$ and $\A(q^{-1})$ such
that 
$$ 
\rho_1({\bf x}){=}{\bf p}, \rho_1({\bf p}){=}{\bf x}, \rho_1({\bf u}){=}
{\bf u}~ {\rm  and} ~\rho_2({\bf x}){=}{\bf x}, \rho_2({\bf p}){=}{\bf p},
\rho_2 ({\bf u}){=}-{\bf u}^\ast.
$$ 
Because the $\ast$-algebras $\A(q)$ and
$\A(q^{-1})$ are isomorphic, we shall assume in what follows that $0<q<1$.
\bn
{\bf  II. An operator-theoretic model}\mn
{\bf II.1.}~ Let $\mu_1$ be a finite positive Borel measure on the intervall
$[q,1)$. The measure $\mu_1$ extends uniquely to a Borel measure $\mu$ on the 
half-axis $\bbbr_+=(0,+\infty)$ by setting  $\mu (q^n\M):=q^n\mu_1(\M)$ for
any Borel subset $\M$ of $[q,1)$. Then $\mu$ has obviously the property
that $\mu (q\N)=q\mu (\N)$ for an arbitrary Borel subset $\N$ of $\bbbr_+$
or equivalently that ${{d\mu (at)}\over {qt}}={{d\mu(t)}\over{t}}$ for $t\in
\bbbr_+$. We shall work with the Hilbert spaces $\H :=L^2(\bbbr_+,\mu)$ 
and $\gH:=L^2([q,1),\mu_1)$.
First we define three linear operators $U,P$ and $X$ on the Hilbert space
$\H$:\mn 

\litem{(i)} $(Uf)(t)=q^{1/2} f (qt)$ for $f{\bf \in\H}$,
\litem{(ii)} $(Pf)(t)=t f(t)$ for $f{\bf \in\D} (P):=\{f\in\H
 :t f (t)\in\H\},$
\litem{(iii)} $(X f)(t)={\rm i}~t^{-1} ({{f(q^{-1} t)- f(qt)}})$ for
$f \in \D (X) := \{f\in\H :t^{-1} f(t)\in\H\}.$\mn

These operators will play a crucial role throughout this paper. 
Roughly speaking and ignoring technical subtleties (domians, boundary conditions
etc.), 
we shall show that for all "well-behaved" $\ast$-representations of the $q$-defomred 
Heisenberg algebra $\A(q)$ the images of the generators 
${\bf u},{\bf p}$  and ${\bf x}$ act by the same formulas as 
the operators $U, P$ and $X$, respectively. 

Obviously, $P$ is an unbounded self-adjoint operator on $\H$. Using the
relation ${{d\mu (q t)}\over {qt}}={{d\mu (t)}\over t}$ one easily verifies
that $U$ is a unitary operator and that $X$ is a symmetric operator
on $\H$. Let $\D_0$ be the set of functions $f\in\H$ such that supp $f\in
[a,b]$ for some $a>0$ and $b>0$. (Note that $a$ and $b$ may depend on $f$.)
Clearly, $\D_0$ is dense linear subspace of $\H$ which is invariant under
$U,P$ and $X$. It is straightforward to check that 
the operators $P,X,U$ applied to functions $f\in\D_0$ 
satisfy the defining relations (1), (2) and (3)
of the $\ast$-algebra $\A(q)$. In turns out that the symmetric operator $X$ is not
essentially self-adjoint. Our next aim is to characterize the domain of 
the adjoint operator $X^\ast$.

For $f\in \gH = L^2([q,1),\mu_1)$ let $f^e$ and $f^o$ be the functions on
$\bbbr_+$ defined by
$$
f^e(q^{2n} t)=f^o(q^{2n+1}t)=f(t)~{\rm for}~n\in\bbbn_0,~
t\in [q,1)~{\rm and }~f^e(t)=f^o(t)=0~{\rm otherwise}.\eqno (4)
$$ 
Clearly, $f^e$ and $f^o$ are in $\H=L^2
(\bbbr_+,\mu)$ and we have $U(f^e)-q^{1/2} f^0\in\D(X)$ and $Uf^o-q^{1/2}
f^o\in\D(X)$. Let $\gH_e$ and $\gH_o$ denote the set of functions $f^e$
and  $f^o$, respectively, where $f\in\gH=L^2([q,1),\mu_1)$.
\bn
{\bf Lemma 1.} {\it The domain $\D (X^\ast)$ is the direct sum of vector spaces
$\D(X),\gH_e$ and $\gH_o$.}\bn
{\bf Proof.}
It is straightforward to check that $\D(X)+\gH_e+\gH_o\subseteq\D (X^\ast)$.
In order to prove the converse, let $g\in\D(X^\ast)$. Then, by definition
there is an $h\in\H$ such that $\langle Xf,g\rangle=\langle f,h\rangle$
for all $f\in\D(X)$. Inserting the definition of $X$ and using once more the
fact that ${ {\rm {d\mu (qt)} \over {qt} }}={{d\mu (t)}\over t}$ we easily conclude
that $h(t)= = {\rm i} t^{-1}(g(q^{-1} t)-g(qt))$. For a function $f\in\H$ let
$f_n$ denote the function in $L^2([q,1),\mu^+_1)$ given by $f_n(t)=f(q^n t)$. 
Then we get
$$\| h\|^2_{L^2(\bbbr_+,\mu)}={\sum\limits^\infty_{n=-\infty}}
\| h_n\|^2 q^n\ge {\sum\limits^\infty_{n=0}} { {\|
g_{n+1}-g_{n-1}\|^2} \over {q^{2n}}}~q^n
$$
For $n\in\bbbn$ we set $\alpha_n:=\| g_{n+1}-g_{n-1}\|
q^{-{n\over 2}}$. Since $h\in L^2(\bbbr,\mu)$, the sequence $(\alpha_n)$
is in $l_2$. From the inequality
$$
\| g_{2r} - g_{2s}\|\le \alpha_{2r+1} ~
q^{ {{ {2r+1}\over {2}}}} +\cdots + \alpha_{2s+1}~ q^{{{2s+1}\over 2}}
$$
we obtain
$$
\| g_{2r} - g_{2s}\|^2\le \left({\sum\limits^\infty_{i=2s+1}}
|\alpha_i|^2\right)~q^{2s+1} (1-q^2)^{-1},r\ge s~.\eqno (5)
$$
Since $(\alpha_n)\in l_2$, this implies that the sequence $(g_{2n})_{n\in\bbbn}$ 
converges in the Hilbert space $L^2([q,1),\mu_1)$. Let us denote its limit
by $\xi$. We extend $\xi$ to a function $\xi^e$ on $\bbbr_+$ by setting
$\xi^e(q^{2n}t):=\xi (t)$ and $\xi^e(q^{2n+1} t):=0$ for $n\in\bbbn_0$ ,
$t\in[q,1)$ and $\xi^e(t)=0$ for $t\ge 1$. Replacing even indices by odd
indices, a similar reasoning yields functions $\zeta \in L^2([q,1), \mu_1)$
and $\zeta^o$ on $\bbbr_+$ such that $\zeta^o (q^{2n+1}t)=\zeta (t)$ and
$\zeta^o (q^{2n} t)=0$ for $n\in\bbbn$, $t\in[q,1)$ and $\zeta^o
(t)=0$
for $t\ge 1$. By construction, $\xi^e\in\gH_e$ and $\zeta^o\in\gH_o$.
Our proof is complete once we have shown that $f:=g-\xi^e-\zeta^o$ belongs
to the domain $\D(X)$ of the operator $X$.
 
Letting $r\rightarrow \infty$ in (5), we get
$$
\|\xi-g_{2s}\|^2\le q^{2s+1} (1-q^2)^{-1} {\sum\limits^\infty_
{n=0}} |\alpha_{2n}|^2~.\eqno (6)
$$
From (6) and the corresponding estimation of 
$\|\zeta-g_{2s+1}\|^2$ we obtain 
$$
\eqalign{
&{\sum\limits^\infty_{n=0}}~\| t^{-1}f_n(t)
 \|^2
q^n \le {\sum\limits^\infty_{n=0}}~ {{\|f_n\|^2}\over {q^{2n+2}}}
~q^n
={\sum\limits^\infty_{\tau=0}} ~{{\| \xi-g_{2r}\|^2}\over 
{q^{2r+2}}} + {{\| \zeta-g_{2r+1}\|^2}\over {q^{2r+3}}}\cr
&=(q-q^3)^{-1} {\sum\limits^\infty_{n=0}}~ |\alpha_n|^2 < \infty.\cr}
$$
Since $f (t)=g(t)$ for $t\ge 1$, this inequality implies that the functions
$t^{-1}f(t)$ and $f(t)$ are in $L^2 (\bbbr^+,\mu)$. Thus, $f\in\D
(X)$. $\hfill\qed$
\mn
As shown in the preceding proof, for any function $g\in\D(X^\ast)$ the
"even components" $g_{2n}$ and the "odd components" $g_{2n+1}$ both have 
"boundary limits" $\xi$ and $\zeta$ in $L^2 ([q,1), \mu_1)$. By Lemma 1,
any element $f\in\D(X^\ast)$ is of the form $f=f_X+f^e+f^o$ with uniquely
determined functions $f_X\in \D(X), f^e\in\H_e$ and $f^o\in\H_o$. By the 
definition of $\H_e$ and $\H_o$, there exist unique functions $f_e,f_o\in\gH=
L^2([q,1),\mu_1)$ such that $(f_e)^e=f^e$ and $(f_o)^o=f^o,$ where the
function $(f_e)^e$ and $(f_o)^o$ on $\bbbr$ are given by (4). 
This notation will be kept in the sequel.

Let $\langle\cdot,\cdot\rangle$ and $(\cdot,\cdot)$ denote the scalar products of 
the Hilbert spaces $L^2(\bbbr_+,\mu)$ and $L^2([q,1), t^{-1}\mu_1)$, respectively.
\bn
{\bf Lemma 2.} {\it For arbitrary functions $f,g\in\D (X^\ast)$ we have }
$$
\langle X^\ast f,g\rangle -\langle f,X^\ast g\rangle ={1\over{2{\rm i}}}
\{(f_e+f_o, g_e+g_o)-(f_e-f_o, g_e-g_o)\}.\eqno (7)_+
$$
{\bf Proof.} Let $h\in L^2([q,1),\mu_1)$. From the definitions 
of the operator $X$ and of the functions ${h^e,h^o\in L^2(\bbbr,\mu)}$ 
we easily derive that $(X^\ast h^e)
(t){=}-i t^{-1} h(qt)$ for $t\in[1,q^{-1}), 
(X^\ast h^e)(t){=}0$\break for $t\in\bbbr_+
\backslash [1,q^{-1}),(X^\ast h^o)(t)=-i~ t^{-1} h(t)$ for $t\in[q,1)$ and
$(X^\ast h^o)(t)=0$ for $t\in\bbbr_+\backslash [q,1)$. Inserting these
expressions and using the symmetry of the operator $X$ we compute
$$\eqalign{
&{\ }\qquad\hskip1cm\langle X^\ast f, g\rangle - \langle f, X^\ast g\rangle =
 \langle X^\ast f_o,g_e\rangle -
\langle f_e, X^\ast g_o\rangle\cr
&{\ }\qquad\hskip1cm= -i{\int\limits^1_q} (f_o(t)\overline{g_e(t)} + f_e(t) \overline{g_o(t)})
t^{-1} d\mu (t)\cr
&{\ }\qquad\hskip1cm=-i\{(f_o,g_e)+(f_e,g_o)\}\cr
&{\ }\qquad\hskip1cm={1\over {2{\rm i}}}\{(f_e+f_o, g_e+g_o)-(f_e-f_o, g_e-g_o)\}.\hskip4.5cm\qed\cr}
$$
Let us illustrate the preceding by the simplest example.
\bn
{\bf Example 1.} Let $\mu_1$ be the Delta measure $\delta_a$, where $a$ is a 
fixed number from the intervall $[q,1)$. Then the measure $\mu$ is supported 
on the points $aq^n,n\in\bbbzz$, and we have $\mu (\{a q^n\})=q^n \mu (\{a\})=
q^n$. Hence the scalar product of the Hilbert space $\H=L^2(\bbbr_+,\mu)$ is 
given by the Jackson integral
$$
\langle f,g\rangle = {\sum\limits^{+\infty}_{n=-\infty}} f(aq^n) ~\overline{
g(aq^n)} ~q^n.
$$
Let $e_n\in\H$ be the function $e_n(t)=q^{ { {-n}\over 2} } \delta^t_{aq^n}$, where 
$\delta^t_s$ is the usual Kronecker symbol. Then the 
vectors $e_n, n\in\bbbzz$, form an orthonormal basis of $\H$ and the
actions of the operators $U,P,X$ on these vectors are given by
$$
U e_n = e_{n-1}, ~P e_n=aq^n e_n,~ X e_n= {{\rm i}\over{aq^n}}\left( q^{-1/2} 
e_{n+1} -q^{1/2} e_{n-1}\right)~.
$$
These equations are in accordance with formulas (5) in ref. 3. 
If $f$ is the function
in $L^2([q,1),\mu_1){\cong}\bbbc$ with $f(a){=}1$, then by definition
$f^e(aq^{2n})=f^o(aq^{2n+1})=1,f^e(aq^{2n+1})=$\break
$f^o(aq^{2n})=0$ for $n\in\bbbn_0$
and $f^e(t)=f^o(t)=0$ for $t\ge 1$. Then we have $\D(X^\ast)=\D(X)+\bbbc
\cdot f^e+\bbbc\cdot f^o$ by Lemma 1 and formula (7)$_+$ reads as
$$\eqalign{
&\langle X^\ast (\varphi+\alpha_1f^e+\beta_1 f^o),
\psi+\alpha_2f^e+\beta_2 f^o
\rangle-\langle\varphi+\alpha_1 f^e+\beta_1f^o,
X^\ast(\psi+\alpha_2 f^e+\beta_2 
f^o)\rangle\cr
&=-{\rm i}~ a^{-1}\{ \beta_1 \bar{\alpha_2} + \alpha_1 {\bar{\beta_2}} \}
= {1\over {2{\rm i}a}} \{ (\alpha_1 + \beta_1) \overline{(\alpha_2+\beta_2)}
-(\alpha_1-\beta_1)\overline{(\alpha_2-\beta_2)}\}\cr}
$$
for $\alpha_1,\beta_1,\alpha_2,\beta_2\in\bbbc. $ \hfill\qed
\bn
{\bf II.2} The above considerations carry over almost verbatim to the
case where the positive half-axis $\bbbr_+$ is replaced by the negative
half-axis $\bbbr_-=(-\infty,0)$. Any positive finite Borel measure $\mu_1$ 
on the intervall $[q,1)$ induces a positive Borel measure $\mu$ on $\bbbr_-$
by defining $\mu(-q^n \M):=q^n\mu_1(\M)$ for a Borel subset $\M$ of $[q,1)$.
The operators $U,P,X$ on the Hilbert space $\H_-:=L^2(\bbbr_-,\mu)$ are
defined by the same formulas as in the preceding subsection and Lemma 1 and
 its proof remain valid in this case as well. However, there is an essential 
difference which will be crucial in the sequel: Since in the proof of Lemma 2
the integration is over the intervall $(-1,-q]$, the expression on the
right hand side of $(7)_+$ must be multiplied by $-1$. That is, instead
of $(7)_+$ we now have
$$
\langle X^\ast f,g\rangle -\langle f,X^\ast g\rangle = 
{1\over {2{\rm i}}} \{ (f_e+f_o,g_e+g_o)-(f_e-f_o,g_e-g_o)\}\eqno (7)_-
$$
for $f,g\in\D (X^\ast).$
\bn
{\bf II.3} After the preceding preparations we are now able to develop the 
operator-theoretic model for the description of $\ast$-representations
of the $q$-Heisenberg algebra $\A(q)$. For this let us fix two families
$\{\mu^{j,+}_1; j\in I_+\}$ and $\{\mu^{j,-}_1; j\in I_-\}$ of finite 
positive Borel measures on the intervall $[q,1)$.\nz
As above, we define the Hilbert spaces $\H_{j,\pm} :=L^2 (\bbbr_\pm,\mu^{j,\pm}),
~j\in I_\pm$, and the operators $U_{j,\pm}, P_{j,\pm},X_{j,\pm}$ acting
therein. We shall work with the representation Hilbert space $\H=\H_+
\oplus\H_-$, where $\H_+:={\bigoplus\limits_{j\in I_+}} \H_{j,+}$ and $\H_-
:={\bigoplus\limits_{j\in I_-}}\H_{j,-}$. The elements of $\H$ are pairs
$\gf=(\gf^+,\gf^-)$, where $\gf^+=(f^{j,+};j\in I_+)\in\H_+$ and 
$\gf^-=(f^{j,-};j\in I)\in \H_-$. Let $U,P,X$ denote the operators on 
$\H$ which are
defined as the direct sums of the operators $U_{j,+}, U_{j,-}; P_{j,+},
P_{j,-}; X_{j,+}, X_{j,-}$, respectively. Clearly, $U$ is a unitary 
operator and $P$ is a self-adjoint operator on $\H$. The operator $X$ is
only symmetric, but not self-adjoint. Our next aim is to describe all
self-adjoint extensions $\tilde{X}$ of $X$ on $\H$ which have the property 
that $U\tilde{X}U^{-1}=q\tilde{X}$. 

Let $V$ and $W$ be two unitary linear transformations of
the Hilbert space $\gH_- :={\bigoplus\nolimits_{j\in I_-}}
L^2 ([q,1), t^{-1}\mu^{j,-}_1)$ on the Hilbert space 
$\gH_+ :={\bigoplus\nolimits_{j\in I_+}} L^2 ([q,1), t^{-1} \mu^{j,+}_1)$
. We define a linear operator $X_{V,W}$ as
being the restriction of the adjoint operator $X^\ast$ to the domain
$$\eqalign{
&\D(X_{V,W}):= \{ \gf=\gf_X + \gf^e + \gf^o \in \D (X^\ast) : \gf_X \in \D (X),
\cr
&\gf^+_e = V (\gf^-_e + \gf^-_0) + W (\gf^-_e - \gf^-_0),\, 
\gf^+_0 = V (\gf^+_e + \gf^-_0) - W(\gf^-_e - \gf^-_0).\cr}\eqno(8)
$$
\bn
{\bf Proposition 3.} {\it $X_{V,W}$ is a self-adjoint operator on $\H$ such that
$X{\subseteq} X_{V,W}$ and $UX_{V,W}U^\ast{=}$\break $q~ X_{V,W}$. In 
particular, we have $U\D(X_{V,W})=\D(X_{V,W})$. Conversely, for any 
self-adjoint extension $\tilde{X}$ of
$X$ satisfying $U\D(\tilde{\X})\subseteq\D (\tilde{X})$ there exist unitary
transformations $V,W$ of $\gH_+$ onto $\gH_-$ such that $\tilde{X}=X_{V,W}$.}
\bn
{\bf Proof.} From (7)$_+$ and (7)$_-$ we obtain
$$\eqalign{
&-2i (\langle X^\ast \gf,\gg\rangle - \langle \gf,X^\ast \gg\rangle)\cr
&={(\gf^+_e+\gf^+_o,\gg^+_e+\gg^+_o)+(\gf^-_e-\gf^-_o,\gg^-_e-\gg^-_o)
-(\gf^-_e+\gf^-_o,\gg^-_e+\gg^-_o)-(\gf^+_e-\gf^+_o,\gg^+_e-\gg^+_o)}.\cr}\eqno (9)
$$
for arbitrary elements $\gf=\gf_X+\gf^e+\gf^o$ and $\gg=\gg_X+\gg^e+\gg^o$ of 
$\D(X^\ast)$.
Here $\gf^+_e$ denotes the sequence $(f^{j,+}_e;j\in I_+)\in\gH_+$ with $
f^{j,+}_e\in L^2([q,1),\mu^j_1)$ such that the extension $(f^{j,+}_e)^e$ of
$f^{j,+}_e$ to $\bbbr_+$ by means of formula  (4) is just the $(j,+)$-component
of the vector $\gf^e\in\H$. A similar meaning attached to the other symbols
$\gf^-_e,\gf^+_o,\gf^-_o,\gg^+_e,\gg^-_e,\gg^+_o,\gg^-_o$ occuring in (9).
If $\gf,\gg\in\D (X_{V,W})$, then we have $\gf^+_e+\gf^-_o=V(\gf^-_e+\gf^-_o)$,
$\gg^+_e+\gg^-_o=V(\gg^-_e+\gg^-_o),\gf^+_e-\gf^+_o=W(\gf^-_e-\gf^-_o)$ and 
$\gg^+_e-\gg^+_o=W(\gg^-_e-\gg^-_o)$ by (8). Since $X_{V,W}\subseteq X^\ast,$ 
we therefore obtain that $\langle X_{V,W}\gf,\gg\rangle-\langle \gf,X_{V,W}\gg
\rangle =0$ by (9), that is, the operator $X_{V,W}$ is symmetric. 
Now let $\gg\in\D((X_{V,W})^\ast)$. Since $X\subseteq X_{V,W}\subseteq (X_{V,W})^\ast
\subseteq X^\ast$, we then have $\langle X^\ast\gf,\gg\rangle=\langle 
\gf,X^\ast\gg\rangle$ and hence
$$
(\gf^+_e+\gf^+_o,\gg^+_e+\gg^+_o)+(\gf^-_e-\gf^-_o,\gg^-_e-\gg^-_o)
=(\gf^-_e+\gf^-_o,\gg^-_e+\gg^-_o)+(\gf^+_e-\gf^+_o,\gg^+_e-\gg^+_o)\eqno (10)
$$
for all $\gf\in\D(X_{V,W})$ by (9). Inserting (8) into (10), we get 
$$\eqalign{
&(f^-_e+\gf^-_o,V^\ast(\gg^+_e+\gg^-_o))+(\gf^-_e-\gf^-_o,\gg^-_e-\gg^-_o)\cr
&=(\gf^-_e+\gf^-_o,\gg^-_e+\gg^-_o)+(\gf^-_e-\gf^-_o,W^\ast(\gg^+_e-\gg^+_o)).\cr}\eqno (11)
$$
From the construction it is clear that for arbitrary $\h,\k\in\gH_-$ there 
exists $\gf\in\D(X_{V,W})$ such that $\gf^-_e+\gf^-_o=\h$ and $\gf^-_e-\gf^-_o=\k$. 
Therefore, it follows from (11) that $V^\ast(\gg^+_e+\gg^+_o)=\gg^+_e+\gg^-_o$ and 
$W^\ast(\gg^+_e-\gg^+_o)=\gg^-_e-\gg^-_o$ which in turn implies that $\gg\in\D(X_{V,W})$.
Thus we have shown that the operator $X_{V,W}$  is self-adjoint. From the
relations $U(\gf_e)-q^{1/2}\gf_o\in\D(X)$ and $U(\gf_o)-q^{1/2}\gf_e\in\partial
(X)$ we see that $U\D(X_{V,W})=\D(X_{V,W})$. Since $UX U^\ast=qX$ and hence
$UX^\ast U^\ast=qX^\ast$ and $X_{V,W}$ is the restriction of $X^\ast$ to
$\D (X_{V,W})$, the latter yields $UX_{V,W} U^\ast=qX_{V,W}$.

Conversely, suppose that $\tilde{X}$ is a self-adjoint extension of $X$ such
that $U\D(\tilde{X})\subseteq\D(\tilde{X})$. Since $\tilde{X}$ is symmetric,
we have equation (10) for arbitrary elements $\gf,\gg\in D(\tilde{X}).$
By assumption, $U\gf\in\D(\tilde{X})$ for all $\gf\in\D(\tilde{X}).$ 
Replacing $\gf$ by $U\gf$ in (10) we get
$$\eqalign{
&(\gf^+_e+\gf^+_o,\gg^+_e+\gg^+_o)+(\gf^-_o-\gf^-_e,\gg^-_e-\gg^-_o)\cr
&=(\gf^-_e+\gf^-_o,\gg^-_e+\gg^-_o)+(\gf^+_o-\gf^+_e,\gg^+_e-\gg^+_o).\cr}\eqno (12)
$$
Setting $\gf=\gg$ and combining formulas (10) and (12) we obtain
$$
\| \gf^+_e+\gf^+_o\| =\| \gf^-_e+\gf^-_o\|~{\rm and}~\|\gf^+_e-\gf^+_o\| =\|
\gf^-_e -\gf^-_o\| \eqno (13)
$$
for all $\gf\in\D(\tilde{X})$.

For $\gf\in\D(X^\ast)$ we abbreviate $B_\pm(\gf)=(\gf^\pm_e+\gf^\pm_o,
\gf^\pm_e-\gf^\pm_o)$. The vector space $B_\pm(\tilde{X})=\{ B_\pm (\gf) :\gf\in
\D(\tilde{X})\}$ is called the "boundary space" of the operator $\tilde{X}$.
We shall show that $B_+(\tilde{X})=\gH_+\oplus\gH_+$ and $B_-(\tilde{X})=
\gH_-\oplus\gH_-$. First let us note that the spaces $B_\pm (\tilde{X})$ are
closed in $\gH_\pm\oplus\gH_\pm$. Otherwise let $\tilde{\tilde{X}}$ denote
the restriction of $X^\ast$ to the domain $\D(\tilde{\tilde{X}})=\{\gf\in\D
(X^\ast):B_\pm (\gf)\in \overline{B_\pm(\tilde{X})}\, \}$, where the bar means 
the closure in the Hilbert space $\gH_\pm\oplus\gH_\pm$. The symmetry of an 
operator $Y$ such that $X\subseteq Y\subseteq X^\ast$ is equivalent to the
validity of equation (10) for all $\gf,\gg\in\D(Y)$. Hence $\tilde{\tilde{X}}$
is symmetric, because $\tilde{X}$ is so. Since a self-adjoint operator has
no proper symmetric extension, we conclude that 
$\tilde{X}=\tilde{\tilde{X}}$ which means
that $B_+(\tilde{X})$ and $B_-(\tilde{X})$ are closed. Next let us suppose
that $(\xi, \zeta)\bot B_+(\tilde{X})$ in $\gH_+\oplus\gH_+$. We then choose
a vector $\gg\in\D(\tilde{X})$ such that $\xi=\gg^+_e+\gg^+_o, \zeta=
\gg^+_e-\gg^+_o$ and $\gg^-_e=\gg^-_o=0$.
Then the right-hand side of (9) vanishes for all $\gf\in\D(\tilde{X})$, so
that $\langle\tilde{X}\gf,\gg\rangle=\langle X^\ast\gf,\gg\rangle =\langle
\gf,X^\ast\gg\rangle$ for all $\gf\in\D(\tilde{X})$ by (9). Consequently,
$\gg\in\D(\tilde{X}^\ast)$. Since $\tilde{X}$ is self-adjoint, $\gg$ must be in $\D(\tilde{X})$. Because $(\xi,\zeta)\bot B_+(\tilde{X})$, this implies
that $\xi=\zeta=0$. This proves that $B_+(\tilde{X})=\gH_+\oplus\gH_+$. Similarly
$B_-(\tilde{X})=\gH_-\oplus\gH_-$.

Since $B_\pm (\tilde{X})=\gH_\pm\oplus\gH_\pm$ as just shown, is follows
from (13) that there are unitary operators $V$ and $W$ of $\gH_-$ onto
$\gH_+$ such that $\gf^+_e+\gf^+_o=V(\gf^-_e+\gf^-_o)$ and $\gf^+_e-\gf^+_o=
W(\gf^-_e-\gf^-_o)$ for all $\gf\in\D(\tilde{X})$. That is, $\D(\tilde{X})
\subseteq\D(X_{V,W})$. Since $\tilde{X}$ and $X_{V,W}$ are self-adjoint, we 
conclude that $\tilde{X}=X_{V,W}$. $\hfill\qed$
\bn
{\bf III. Irreducibility and unitary equivalence}\mn
{\bf III.1} The next two propositions decide when a triple of operators 
$\{P,X_{V,W},U\}$ defined in the preceding section is irreducible and when
 two such triples are unitarily equivalent. Here we shall say that the
triple $\{P,X_{V,W},U\}$ on $\H$ is {\it irreducible} if any
bounded operator $A$ on $\H$ satisfying 
$$
PA\subseteq AP,~X_{V,W} A\subseteq A X_{V,W} ~{\rm and}~ AU=UA \eqno (14)
$$
is a scalar multiple of the identity operator on $\H$. 

Recall that the operator triple 
$\{P,X_{V,W}, U\}$ depends on the two families $\{\mu^{j,\pm}_1;j\in 
I_\pm\}$ of measures on the intervall $[q,1)$ and on the two unitary operators 
$V,W:\gH_-\rightarrow\gH_+$. In order to formulate the corresponding conditions
it is convenient to work with the Hilbert spaces $\gK_\pm={\bigoplus\nolimits_{j\in I_\pm}}
L^2([q,1),\mu^{j,\pm}_1)$ rather than with $\gH_\pm={\bigoplus\nolimits_{j\in I_\pm}}
L^2([q,1),t^{-1}\mu^{j,\pm}_1)$. Further, let $P_\pm$ denote the self-adjoint
operator on $\gK_\pm$ which acts componentwise as the multiplication by the 
variable $t$. Clearly, $V$ and $W$ are bounded linear operators of $\gK_-$ to $\gK_+$
such that 
$$
V^\prime := P^{1/2}_+~ V P^{-1/2}_-~~{\rm and}~~W^\prime :=P^{1/2}_+ ~WP^{-1/2}_-\eqno (15)
$$
are unitary. 
\mn
{\bf Proposition 4.} {\it The triple $\{P,X_{V,W},U\}$ as defined above is irreducible
if and only if any bounded self-adjoint operators $A_+$ on $\gK_+$ and 
$A_-$ on $\gK_-$ satisfying 
$$
A_+P_+=P_+A_+, A_-P_-=P_-A_-, A_+V^\prime= V^\prime A_-,A_+W^\prime = 
W^\prime A_-\eqno (16)
$$
or equivalently
$$
A_+P_+=P_+A_+,A_-P_-=P_-A_-,A_+ V=VA_+, A_+W=WA_-\eqno (17)
$$
are scalar multiples of the identity.}\mn
{\bf Proof.} We only show that the above condition implies the irreducibility
of the triple. The proof of the converse implication is easier and will
be omitted. Suppose that $A$ is a bounded operator on $\H$ satisfying (14).
Since the set of such $A$ is invariant under the involution, we can assume
that $A$ is self-adjoint. Let $E(\cdot)$ denote the spectral projections
of $P$. Since $PA\subseteq AP$, the subspace $\gK_+=E([q,1))\H$ of $\H$ 
reduces $A$ and the restriction $A_+$ of $A$ to $\gK_+$ commutes with the
restriction $P_+$ of $P$ to $\gK_+$. Similarly, the restrictions $\tilde{A}_-$
of $A$ and $\tilde{P}_-$ of $P$ to the reducing subspace $E((-1,q])\H$
commute. Changing the variable from $t$ to $-t$, the Hilbert space 
$E((-1,q])\H$ and the operator $\tilde{P}_-$ become $\gK_-$ and $P_-$,
respectively, and the operator $\tilde{A}_-$ goes into an operator, say
 $A_-$, on $\gK_-$. Thus, $A_- P_-=P_-A_-$. From the assumptions $AU=UA$
and $X_{V,W}A\subseteq AX_{V,W}$ it follows easily that $(A\gf)^\pm_e=
A_\pm\gf^\pm_e$ and $(A\gf)^\pm_o=A_\pm\gf^\pm_o$ for $\gf\in\D(X_{V,W})$. Since
$A\gf\in\D(X_{V,W})$ has to satisfy the relation (8), we obtain $A_+V=VA_-$ and 
$A_-W=WA_-$. Therefore, by the above condition, $A_\pm=\lambda_\pm I$ for
some $\lambda_\pm\in\bbbc$. Since $A_+V=VA_-$ and $AU=UA$, it follows that
$\lambda_+=\lambda_-$ and $A=\lambda_+\cdot I$ on $\H$. $\hfill\qed$ 

Using similar operator-theoretic arguments it is not difficult to prove
\mn
{\bf Proposition 5.} {\it Two triples $\{P,X_{V,W}, U\}$ and $\{\tilde{P},
\X_{\tilde{V},\tilde{W}},\tilde{U}\}$ are unitarily equivalent
if and only if there unitary operators $A_+$ of $\gK_+$ to $\tilde{\gK}_+$ and
$A_-$ of $\gK_-$ to $\tilde{\gK}_-$ such that
$$
A_+P_+=\tilde{P}_+A_+,A_-P_-=\tilde{P}_-A_-,A_+V=\tilde{V}A_-~{\rm and}~
A_+W=\tilde{W} A_-~,\eqno (18)
$$
where the tilde refers to the corresponding operators and spaces for the
triple \break $\{\tilde{P}, X_{\tilde{V},\tilde{W}},\tilde{U}\}$.}
\bn
{\bf III.2.} We shall illustrate the preceding by describing a few examples of
irreducible representations. We begin with the simplest possible case.\mn
{\bf Example 2.} Suppose that the Hilbert spaces $\gK_+$ and $\gK_-$ are
one-dimensional. Then the families of measure $\{\mu^{j,+}_i;\in I_+\}$ and 
$\{\mu_1^{j,-}; j\in I_-\}$ consist only of single Dirac measures $\delta_a$
and $\delta_b$, respectively, where $a,b\in[q,1)$. Then the triples $\{
P,X_{V,W},U\}$ are parametrized by complex numbers $V=V^\prime=e^{i\varphi}$
and $W=W^\prime=e^{i\psi},\varphi,\psi\in\bbbr$. The 
self-adjoint extension $X_{V,W}$ is then characterized by the boundary
condition (8), that is,
$$
\gf^+_e+\gf^+_o=e^{i\varphi}(\gf^-_e+\gf^-_o),~ \gf^+_e-\gf^+_o=
e^{i\psi}(\gf^-_e-\gf^-_o).
$$
Each such triple is irreducible because the condition in Proposition 4 is
trivially fulfilled. Two triples with different pairs of numbers $(V,W)$ are 
not unitary equivalent. The case where $e^{i\varphi}=e^{i\psi}=1$ and $a=b$
has been treated in detail in ref. 3. \hfill\qed
\mn
{\bf Example 3.} Let $P_+$ be a self-adjoint operator and $Z$ a unitary
operator on a Hilbert space $\gK_+$ such that the commutant $\{P_+,Z\}^\prime$
is equal to $\bbbc\cdot I$. Such operators exist 
on any separable Hilbert space$^{12}$. Upon scaling we can assume that 
the spectrum of $P_+$ is contained in $[q,1)$. By the
spectral representation theorem$^{1,ch. X, 5.}$, we can represent $P_+$ up to
unitary equivalence as the multiplication operator by the independent variable
$t$ on some direct sum Hilbert space $\gK_+={\bigoplus\nolimits_{j\in I_+}} L^2
([q,1);\mu^{j,+}_1)$. Let $\{\mu^{j,-}_1; j\in I\}$ be an arbitrary family
of measures on $[q,1)$ such that $\dim\gK_+=\dim\gK_-,$ where $\gK_-:=
{\bigoplus\nolimits_{j\in I_-}} L^2([q,1);\mu^{j,-}_1)$. Let $W^\prime$ be a
unitary operator from $\gK_-$ to $\gK_+$. We set $V^\prime:= ZW^\prime$ and define
$V$ and $W$ by (15). Then the triple $\{P,X_{V,W},U\}$ is irreducible. 

Indeed, if $A_+$ and $A_-$
be bounded self-adjoint operators satisfying (17), then we have $A_+Z=A_+
V^\prime W^{\prime\ast}=V^\prime A_-W^{\prime\ast}=V^\prime W^{\prime\ast}
A_+=ZA_+$ and $A_+P_+=P_+A_+$, so that $A_+=\lambda \cdot I$ for some $\lambda
\in\bbbc$ and hence $A_-=V^{\prime\ast} A_+V^\prime=\lambda\cdot I$.
By Proposition 4, the triple is irreducible.\hfill\qed
\mn
{\bf Example 4.} For this example we assume that there exist numbers $a,b\in
[q,1)$ such that $\mu^{j,+}_1=\delta_a$ and $\mu^{k,-}_1=\delta_b$ for all
$j\in I_+$ and $k\in I_-.$ We shall show that in this case an irreducible
triple $\{P,X_{V,W}, U\}$ can be only obtained if both index sets $I_+$ and
$I_-$ are singletons or equivalently if $\dim \gK_+=\dim \gK_-=1$.
Indeed, otherwise we take a self-adjoint operator $A_+$ on $\gK_+$ such
that $A_+ V^\prime W^{\prime\ast}=V^\prime W^{\prime\ast} A_+$ and $A_+
\not\in\bbbc\cdot I$ and set $A_-:=V^{\prime\ast}A_+V^\prime$. Then the
conditions (16) are fulfilled, hence the triple is not irreducible. 
\hfill\qed
\mn
{\bf Example 5.} If the spectra of the operators $P_+$ on $\gK_+$ and
$P_-$ on $\gK_-$ are singletons, then we have seen in Example 4 that 
irreducible triples exist only in the trivial case where $I_+$ and $I_-$
are singletons. We now show that this is no longer true if both spectra consist
of two points. To be more precise, we shall consider the following situation:
The index sets $I_\pm$ are disjoint union of two countable infinite sets
$I^1_\pm$ and $I^2_\pm$ and there are numbers $a_1,a_2,b_1,b_2\in [q,1)$,
$a_1\ne a_2$, such that $\mu^{j,+}_1 =\delta_{a_1}$ for $j\in I^1_+,
\mu^{j,+}_1=\delta_{a_2}$ for $j\in I^2_+,\mu^{j,-}_1=\delta_{b_1}$ for $j\in I^1_-$
and $\mu^{j,-}_1=\delta_{b_2}$ for $j\in I^2_-$. By identifying $I^j_\pm$
with the natural numbers the Hilbert spaces $\gK_+$ and $\gK_-$ become
the direct sum $l_2(\bbbn)\oplus l_2(\bbbn)$ of two $l_2$-spaces. We choose
a bounded operator $T$ on $l_2(\bbbn)$ such that $\{T,T^\ast\}^\prime=
\bbbc\cdot I$ and $I\le 3 T^\ast T\le 2\cdot I$ and $I\le 3T^\ast T\le 2\cdot
I$. It is well-known (see ref. 8, Anhang, \S 4) that the operator matrix
$$
Z=\left( \matrix{ T \qquad\quad \sqrt{I-TT^\ast}\cr
-\sqrt{I-T^\ast T}\qquad\quad T^\ast\cr}\right)
$$
defines a unitary operator $Z$ on $\gK_+=\gK_-=l_2(\bbbn)\oplus l_2(\bbbn).$
Let $W^\prime$ be an arbitrary unitary operator on $\gK_+=\gK_-$ and set $V^\prime
:=ZW^\prime$. Then the triple $(P,X_{V,W},U)$ is irreducible.

Indeed, let $A_+$ and $A_-$ be self-adjoint bounded operators on $\gK_+=\gK_-$
satisfying (17). Since $a_1\ne a_2$, the relation $A_+P_+=P_+A_+$ implies
that $A_+$ is given by a diagonal operator matrix
$$
A_+=\left(\matrix{ B&0\cr 0&C\cr}\right).
$$
From (17) we get $A_+Z=ZA_+$. Comparing the matrix entries of the first line yields
$BT=TB$ and $B\sqrt{ I-T^\ast T}=\sqrt{I-T^\ast T}~C$. Since $B=B^\ast$, we
have $BT^\ast=T^\ast B$. Therefore, $B$ commutes with $T$ and $T^\ast$ and
so with $\sqrt{I-T^\ast T}$ which in turn gives $\sqrt{I-T^\ast T}~B=\sqrt{I-T^\ast T}~ C$.
Because $\sqrt{I-T^\ast T}$ is invertible, we get $B=C$. Since $B\in\{T,T^\ast\}^\prime$
, we obtain $B=C=\lambda\cdot I$ for some $\lambda\in\bbbc$. Thus, $A_+=\lambda
\cdot I$ and $A_-=V^{\prime\ast} A_+V=\lambda\cdot I$, so that the
triple is irreducible by Proposition 4.\hfill\qed
\bn
{\bf IV. A characterization of the operator triples}\mn
Let $\{P,X_{V,W}, U\}$ be an operator triple as in section II and let
$\D_1$ be the set of all vectors $\gf=\gf_X+\gf^e+\gf^o\in\D (X_{V,W})$ with
$\gf_X\in\D_o$, where $\D_o$ is as defined in Section II. Then $\D_1$
is a dense linear subspace of the Hilbert space $\H$ such that $\D_1$
is invariant under the operators $P$, $X_{V,W},U$ and the restrictions of $P$ and 
$X_{V,W}$ to $\D_1$ are essentially self-adjoint. Further, the three 
operators $P,X_{V,W}, U$  applied to vectors $\gf\in\D_1$ satisfy the 
relations (1) and (2). From the construction it is clear that the range 
$E([q,1))\H(\cong\gK_+)$ of the spectral projection $E([q,1))$ of the
operator $P$ is contained in $\D_1$. Our
next proposition says that the operator triples $\{P,X_{V,W},U\}$ can be
characterized by some of the properties just mentioned.
\bn
{\bf Proposition 6.} {\it Let $\{P^\prime, X^\prime, U^\prime\}$ be a triple of
two self-adjoint operators $P^\prime$ and $X^\prime$ and a 
unitary operator $U^\prime$ on
a Hilbert space $\tilde{\H}$. Let $E(.)$ denote the spectral maesure of
$P^\prime$. Suppose that there exists a linear subspace $\D_1
\subseteq \D(P^\prime X^\prime)\cap\D (X^\prime P^\prime)$ of $\H$ 
such that:\sn
\litem{(i)} $E([q,1))\H\subseteq\D_1$ and $E((-1,-q])\H\subseteq\D_1$. 
\litem{(ii)} The operators $P^\prime,X^\prime,U^\prime$ satisfy the relations (1) and (2) for 
vectors in $\D_1$ .
\litem{(iii)} The restrictions $P^\prime\lceil\D_1$ and $X^\prime\lceil\D_1$
of $P^\prime$ and $X^\prime$ to $\D_1$ are essentially self-adjoint.\sn
Then $\{P^\prime,X^\prime,U^\prime\}$ is unitarily equivalent to an operator
triple $\{P,X_{V,W}, U\}$ defined in Section II.}\mn
{\bf Sketch of proof.} The restriction $P^\prime_1$ of $P^\prime$ to the 
invariant subspace $\H_1:=E([q,1))\tilde{\H}$ is obviously a bounded 
self-adjoint operator on the Hilbert space $\H_1$ with spectrum contained
in the intervall $[q,1]$. By the spectral representation theorem$^1$, there is a family $\{\mu^{j,+}_1; j\in I_+\}$ of finite positive
Borel measures on $[q,1]$ and a unitary isomorphism of $\H_1$ on $\gK_+ :={
\bigoplus\nolimits_{j}} L^2([q,1],\mu_1^{j,+})$ such that $P^\prime_1$ is 
unitarily equivalent to the operator $P_1$ on $\gK_+$ which acts componentwise
as the multiplication by the varia\-ble $t$. Since 1 is not an eigenvalue of
$P^\prime_1$ by construction, we have $\mu^{j,+}_1 (\{1\})=0$ for all $j\in I_+$. For 
simplicity let us identify $\H_1$ with $\gK_+$ and $P^\prime_1$ with $P_1$.

Next we show that $\ker P^\prime=\{0\}$. Let $\gf\in\ker P^\prime$.
Since $P^\prime\lceil\D_1$ is essentially self-adjoint by (iii), there exists a
sequence $\{\gf_n\}$ of vectors $\gf_n\in \D_1$ such that $\gf_n\rightarrow
\gf$ and $P^\prime \gf_n\rightarrow P^\prime \gf=0$ in $\H$. Since $X^\prime P^\prime
\gf_n={\rm i}(q^{1/2} U^{\prime\ast} + q^{1/2} U^\prime)\gf_n$ by (ii) and the
operators $U^\prime$ and $U^{\prime\ast}$ are bounded, we obtain $(q^{-1/2}
U^{\prime\ast}+q^{1/2} U^\prime)\gf=0$ in the limit. This in turn yields
that $q\|\gf\|=\|\gf\|$and so $\gf=0$. 

By (ii), we have $U^\prime P^\prime \gf=qP^\prime U^\prime \gf$ for all 
$\gf\in\D_1$. Since $p^\prime\lceil\D_1$ is essentially self-adjoint, this 
remains valid for $\gf\in\D (P^\prime)$, so that $P^\prime\subseteq 
qU^{\prime\ast}P^\prime U^\prime$. Since $P^\prime$ is self-adjoint, we 
conclude that $P^\prime =qU^{\prime\ast}P^\prime U^\prime$.
Hence we have $U^{\prime n} E(\gN)=E(q^{-n}\gN)$ for any Borel subset $\gN$
of $\bbbr$ and arbitrary $n\in\bbbzz$. Let
$\mu^{j,+}$ be the extension of the measure $\mu^{j,+}_1$ to $\bbbr_+$ as in
II.1. From the preceding considerations it follows that  $E(\bbbr_+)\H=
\oplus_j L^2(\bbbr_+,\mu^{j,+})\equiv\H_+$ and that  $U^\prime$ acts in each
component by formula (i) in subsection II.1. Proceeding in a similar manner, we obtain a 
family $\{\mu^{j,-}_1; j\in I_-\}$ of measures on $[q,1]$ such that $\mu^{j,-}_1
(\{1\})=0$ for $j\in I_-,E(\bbbr_-)\H=\oplus_j L^2(\bbbr_-,\mu^{j,-})\equiv
\H_-$ in the notation of Section II and $U^\prime$ acts componentwise as
  given by formula (i) in II.1. Since $E(\{0\})\H=\ker P^\prime =\{0\}$ as proved in the
preceding paragraph, we conclude that $\H=\H_+\oplus\H_-$. 

From the construction it is clear that $P^\prime$ and $U^\prime$ are the
operators $P$ and $U$, respectively, as in Section II. Let us finally
turn to the operator $X^\prime$. Recall that we have $X^\prime P^\prime \gf=
{\rm i}(q^{-1/2} U^{\prime\ast}+q^{1/2} U^\prime)\gf$ for $\gf\in\D_1$. By arguing as 
the paragraph before last, this relation remains valid for all $\gf\in\D
(P^\prime)$. If $f$ denotes a component of the vector $\gf$, then the
preceding equation yields that $g:=tf\in\H,t^{-1}g=f\in\H$ and 
$(X^\prime g)(t)= 
{\rm i}(q^{-1}f(q^{-1}t)-qf(qt))={\rm i}t^{-1}(g(q^{-1} t)-g(qt))=(Xg)(t)$. Hence
$X^\prime \gf=X\gf$ for all $\gf\in\D (P^\prime)$. Since $X^\prime\lceil
\D_1$ is essentially self-adjoint, the relation $U^\prime 
X^\prime\gf=q^{-1}X^\prime U^\prime\gf$ for $\gf\in\D_1$ by (ii) extends
to vectors $\gf\in\D (X^\prime),$ so that $U^\prime X^\prime U^{\prime\ast}=
q^{-1} X^\prime$.
Thus, $X^{\prime}$ is a self-adjoint extension of the operator $X$ such that 
$U\D(X^\prime)=\D(X^\prime)$. By Proposition 3, $X^\prime$ is 
of the form $X_{V,W}$.
\hfill \qed
\bn
{\bf V. ${\bf \ast}$-Representations of the $q$-Heisenberg algebra}\mn
{\bf V.1} We have considered so far only operator triples and operator relations rather 
than representations of the algebra $\A(q)$. But any operator triple 
$\{P,X_{V,W}, U\}$ gives rise to a self-adjoint representation of the 
$\ast$-algebra as follows. Indeed, let $\D_1$ be the domain defined at the 
beginning of section IV. For vectors in $\D_1$ the operators $P,X_{V,W},U$ 
satisfy the defining relations (1) and (2) of the algebra $\A(q)$. Hence 
there exists a unique $\ast$-representation $\pi_1$ of the $\ast$-algebra 
$\A(q)$ on the domain $\D_1$ such that
$$
\pi_1(p)=P\lceil\D_1, \pi_1(x)=X_{V,W}\lceil\D_1,\pi_1(u)=U\lceil\D_1.
$$
(For the notions on unbounded $\ast$-representations used in what follows we
refer to the monograph$^9$. Recall that the symbol $T\lceil\D_1$ means the restriction 
of $T$ to $\D_1$.)\nz
The $\ast$-representation $\pi_1$ is not yet self-adjoint (see ref. 9, Definition 8.1.10),
because, roughly speaking, $\D_1$ is not the largest possible domain.
However, since the operators $\pi_1(p)$ and $\pi_1(x)$ are essentially
self-adjoint, it follows at once from Proposition 8.1.12 (v) in ref. 9 that the
adjoint representation $\pi:=(\pi_1)^\ast$ is self-adjoint. It is not difficult
to verify that the domain $\D$ of the $\ast$-representation $\pi$ is just the 
intersection of domains of all possible products of the operators $P,X_{V,W},U$
(see ref. 9, Proposition 8.1.17). From these facts
it follows that the operator triple
$\{P,X_{V,W}U\}$ is irreducible if and only if the $\ast$-representation
$\pi$ is so and that two triples are unitarily equivalent if and only if the
corresponding $\ast$-representations are so. That is, Propositions 4 and 5 
provide also the conditions for the irreducibility and the unitary 
equivalence of these $\ast$-representations of the $\ast$-algebra
$\A(q)$. \bn
{\bf V.2} Finally, we briefly discuss how operator representations of the
$q$-deformed Heisenberg algebra $\A(q)$ can be constructed by means of the
Schr"dinger representation $\cP :=-{\rm i}{d\over {dt}}$ and $\cQ:=t$ of the
"ordinary" momentum and position operators.

Let us write $q=e^{-\alpha}$ with $\alpha\in\bbbr$. We define three operators
$U,P,X$ on the Hilbert space $\H =L^2(\bbbr)$:
$$
U=e^{iQ}, P=e^{\alpha \cP}, X= 
{\rm i}(q^{-1/2}e^{-{\rm i}\cQ}-q^{1/2} e^{{\rm i}\cQ})
e^{-\alpha\cP}.\eqno (18)
$$
The vector space $\D:={\rm Lin}\{e^{\gamma t-t^2};\gamma\in\bbbc\}$ is a dense
linear subspace of $\H$. Since the operator $e^{\beta P}, \beta\in \bbbr$,
acts as $(e^{\beta P}f)(t)=f(t-\beta {\rm i})$ on functions $f\in\bbbc$ (see, for
instance, ref. 10 for a rigorous proof), the operators $U,P,X$ satisfy the
relations (1) and (2)$^\prime$ and the restrictions of these operators to the
invariant dense domain $\D$ define a $\ast$-respresentation of the 
$\ast$-algebra $\A(q)$. This operator representation (18) appears already
somewhat hidden in ref. 2. Indeed, if we change the variable $t$ to $e^t$, then
the operator triple $\{ U\oplus U, (-P)\oplus P, (-X)\oplus X\}$ on the direct sum
Hilbert space $\H\oplus\H$ is easily seen to be unitarily equivalent to the 
triple in formula (2.2) in ref. 2.

The operator representation (18) is irreducible on $\H$. Obviously, $U$
is unitary and $P$ is self-adjoint. However, an essential disadvantage
of the representation (18) is that the operator $X$ is only symmetric, but
not essentially self-adjoint. The latter can be shown by the argument used
in the proof of Proposition A.2 in ref. 10. The reason for this failure is the
fact the holomorphic function $h(z)=q^{-1/2} e^{{\rm i}z}-q^{1/2} e^{-{\rm i}z}$ admits 
the zero $z_o={\rm i}{\alpha\over 2}$ in the strip $\{z\in\bbbc :0<{\rm I m} z<\alpha\}$.

\bn
{\bf Reference}\mn
\litem{[1]} Dunford, N. and Schwartz, J.: Linear Operators, Part II. Interscience
Publishers, New York, 1963
\sn
\litem{[2]} Fichtm{\"u}ller, M., Lorek, A. and Wess, J.: $q$-Deformed phase space
and its lattice structure. Preprint MPI-PhI/95-109, Munich, 1995
\sn
\litem{[3]} Hebecker, A., Schreckenberg, S., Schwenk, J., Weich, W. 
and Wess, J.: Representations of a $q$-deformed Heisenberg algebra. 
Z. Phys. C. {\bf 64} (1994), 335 - 359
\sn
\litem{[4]} Koelink,  H.T. and  Swartouw, R.F.: On the zeros of the 
Hahn-Exton$q$-Bessel function and associated $q$-Lommel polynomials. 
J. Math. Anal. and Appl. {\bf 186} (1994), 690--710 
\sn
\litem{[5]} Koornwinder, T.H. and Swartouw, R.F.: On $q$-analogues of the
Fourier and Hankel transforms. Trans. Amer. Math. Soc. {\bf 333} (1992),
445 - 461
\litem{[6]} Lorek, A., Weich, W. and Wess, J.: Non-commutative Euclidean
and Minkowski structures. Preprint, MPI-PhI/96-124, Munich, 1996
\sn
\litem{[7]} Pusz, W. and Woronowicz, S.L.: Twisted second quantization.
Reports Math. Phys. {\bf 27} (1989), 231 - 257
\sn
\litem{[8]} Riesz, F. and Sz.-Nagy, B.: Vorlesungen {\"u}ber Funktionalanalysis.
DVW, Berlin, 1956
\sn
\litem{[9]} Schm{\"u}dgen, K.: Unbounded operator algebras and representation 
theory. Birkh"user, Basel, 1990
\sn
\litem{[10]} Schm{\"u}dgen, K.: Integrable operator representations of $\bbbr^2_q,
X_{q,\gamma}$ and $SL_q(2,\bbbr)$. Commun. Math. Phys. {\bf 159} (1994),
217-237
\sn
\litem{[11]} Schwenk, J. and Wess, J.: A $q$-deformed quantum mechanical toy
model. Phys. Letters {\bf B} {\bf 291} (1992), 273 - 277
\sn
\litem{[12]} Topping, D.M.: Lectures on von Neumann algebras. Van Nostrand,
New York, 1971
\sn
\litem{[13]} Weich, W.: The Hilbert space representations for $SO_q(3)$-symmetric
quantum mechanics. Preprint, Munich, 1994, hep-th /9404029
\sn
\litem{[14]} Wess, J. and Zumino, B.: Covariant differential calculus on the
quantum hyperplane. Nucl. Phys. B. Proc. Suppl. {\bf 18 B} (1991), 302-312
\end